\documentclass[11 pt, twoside, reqno]{amsart}
\usepackage{parskip}
\usepackage{amsthm,amsmath,amssymb,oldgerm}
\usepackage[a4paper]{geometry}
\usepackage{mathrsfs}
\usepackage{verbatim}
\usepackage{hyperref,color}

\theoremstyle{plain}
\newtheorem{thm}{Theorem}[section]
\newtheorem*{theorem*}{Main theorem}

\newtheorem{lem}[thm]{Lemma}

\theoremstyle{definition}

\newcounter{EQNR}
\renewcommand{\contentsname}{Contents\\
{\footnotesize\normalfont(A table
of contents should normally not be included)}}

\title[$L^{\infty}$-norm bounds for Siegel--Jacobi cusp forms]{$L^{\infty}$-norm bounds for Siegel--Jacobi cusp forms}
\date{\today}
\subjclass[2010]{11F11, 11F12}

\begin{document}

\author[A.~Aryasomayajula]{Anilatmaja Aryasomayajula}
\address{Indian Institute of Science Education and Research, Tirupati, Karakambadi Road, Mangalam (P.O.) Tirupati-517507, 
Andhra Pradesh, India}
\email{anilatmaja@gmail.com}

\author[J.~Kramer]{J\"urg Kramer}
\address{Institut f\"ur Mathematik, Humboldt-Universit\"at zu Berlin, Unter den Linden 6, D-10099 Berlin, Germany}
\email{kramer@math.hu-berlin.de}

\author[A-M.~von Pippich]{Anna-Maria von Pippich}
\address{Fachbereich Mathematik und Statistik, Universit\"at Konstanz, Universit\"atsstra{\ss}e 10, D-78464 Konstanz, 
Germany}
\email{anna.pippich@uni-konstanz.de}

\begin{abstract}
In this article, we establish explicit and uniform $L^{\infty}$-norm bounds for $L^{2}$-normalized Siegel--Jacobi cusp forms of 
integral weight $k$ and index $m$ for the Siegel modular group $\Gamma_{0}=\mathrm{Sp}_{2g}(\mathbb{Z})$ for arbitrary 
genus $g\geq 1$. Using the generalization of the classical Eichler--Zagier theta decomposition to higher genus, any such 
Siegel--Jacobi cusp form can be written as a finite linear combination of Siegel cusp forms of half-integral weight $k-1/2$ 
multiplied by the higher-dimensional analogues of the classical Jacobi theta functions. By building upon the uniform $L^{\infty}
$-norm bounds on average for Siegel cusp forms established by J.~Kramer and A.~Mandal~\cite{k1} via the associated 
Bergman kernels, we prove that for $k\in\mathbb{Z}_{\geq g+1}$, $m\in\mathbb{Z}_{\geq 1}$, and a given $\epsilon>0$, 
the $L^{\infty}$-norm bound
\begin{displaymath}
\Vert\phi\Vert_{L^{\infty}}=\sup_{(\tau,z)\in\mathbb{H}_{g}\times\mathbb{C}^{g}}\Vert\phi(\tau,z)\Vert_{\mathrm{Pet}}=O_
{\Gamma_{0},\epsilon}\big(k^{(3g^{2}+5g)/8}\,m^{g^{2}+5g/4+\epsilon}\big)
\end{displaymath}
holds for any Siegel--Jacobi cusp form $\phi$ that is $L^{2}$-normalized with respect to the Petersson inner product. These 
estimates provide the first explicit upper bounds in terms of both parameters $k$ and $m$ for arbitrary genus $g$.
\end{abstract}

\maketitle

\section{Introduction}

\subsection{Background and context}
The study of $L^\infty$-norm bounds of automorphic forms is a long-standing problem at the intersection of harmonic analysis 
and number theory. Estimates for these norms, particularly with respect to parameters such as the weight, level, or index, provide 
a deeper understanding of the pointwise behavior and distribution of automorphic forms. For classical cusp forms on the upper 
half-plane, optimal upper bounds with respect to the weight $k$ have been extensively studied starting with the seminal work of 
H.~Iwaniec and P.~Sarnak~\cite{IS95}.

Jacobi forms, systematically introduced by Eichler and Zagier~\cite{EZ85}, generalize classical modular forms by introducing an 
additional spatial variable in $\mathbb{C}$. An interesting extension of the classical $L^{\infty}$-norm problem is to establish bounds 
for Jacobi forms in terms of both the weight $k$ and the index $m$. P.~Anamby and S.~Das~\cite{AD23} investigated the $L^{\infty}
$-norm problem for individual $L^{2}$-normalized Jacobi cusp forms on the full modular group $\mathrm{SL}_{2}(\mathbb{Z})$ by 
analyzing the Fourier coefficients of Hecke eigenforms. In the recent joint work~\cite{AKP25}, the authors established an alternative 
bound utilizing a geometric and analytic method of proof. To be precise, let $J_{k,m}^{\mathrm{cusp}}(\Gamma_{0})$ denote the 
complex vector space of Jacobi cusp forms of weight $k$ and index $m$ for $\Gamma_{0}=\mathrm{SL}_{2}(\mathbb{Z})$. For a 
Jacobi cusp form $\phi\in J_{k,m}^{\mathrm{cusp}}(\Gamma_{0})$, its pointwise Petersson norm is defined by
\begin{displaymath}
\Vert\phi(\tau,z)\Vert_{\mathrm{Pet}}^{2}:=\vert\phi(\tau,z)\vert^{2}\,\mathrm{Im}(\tau)^{k}e^{-4\pi m\mathrm{Im}(z)^{2}/\mathrm{Im}
(\tau)}\qquad(\tau\in\mathbb{H},\,z\in\mathbb{C}).
\end{displaymath}
Assuming that $\phi$ is $L^{2}$-normalized with respect to the Petersson inner product, the main result in~\cite{AKP25} establishes 
the $L^{\infty}$-norm bound
\begin{displaymath}
\Vert\phi\Vert_{L^{\infty}}:=\sup_{(\tau,z)\in\mathbb{H}\times\mathbb{C}}\Vert \phi(\tau,z)\Vert_{\mathrm{Pet}}=O_{\Gamma_{0},\epsilon}
\big(k\,m^{7/8+\epsilon}\big)
\end{displaymath}
for any $k\in\mathbb{Z}_{\geq 5}$, $m\in\mathbb{Z}_{\geq 1}$, and a given $\epsilon>0$, where the implied constant depends only 
on $\Gamma_{0}$ and the choice of $\epsilon>0$. The method of proof presented in~\cite{AKP25} yields this estimate by using the 
classical Eichler--Zagier theta decomposition. More specifically, we decompose the Jacobi cusp form into a finite linear combination 
of modular forms of half-integral weight $k-1/2$ multiplied by classical Jacobi theta functions, linking their $L^{\infty}$-norms to the 
geometry of the Bergman kernel. Although the resulting dependence on the index $m$ slightly differs from the bound $O_{\epsilon}
\big((km)^{1+\epsilon}\big)$ obtained by P.~Anamby and S.~Das~\cite{AD23}, this structural approach has the advantage of naturally
extending to higher genus.

The primary objective of the present paper is to extend this geometric and analytic methodology to the higher-dimensional setting. 
Specifically, we investigate the $L^{\infty}$-norm problem for Siegel--Jacobi cusp forms of weight $k$ and index $m$ for the Siegel 
modular group $\mathrm{Sp}_{2g}(\mathbb{Z})$ for arbitrary genus $g\geq 1$. Extending the case $g=1$, any such Siegel--Jacobi 
form admits an analogous theta decomposition into Siegel modular forms of half-integral weight and higher-dimensional Jacobi theta 
functions~\cite{EZ85, ziegler}. By controlling the components of this decomposition, we establish our explicit bounds for $L^{2}
$-normalized Siegel--Jacobi cusp forms in terms of both the weight $k$ and the index $m$ for arbitrary genus $g\geq1$.

\subsection{Main result}
The main result of this paper provides explicit and uniform $L^{\infty}$-norm bounds for $L^{2}$-normalized Siegel--Jacobi cusp forms 
on the Siegel modular group $\Gamma_{0}=\mathrm{Sp}_{2g}(\mathbb{Z})$ for arbitrary genus $g\geq 1$. To formulate our result
precisely, let $\mathbb{H}_{g}$ denote the Siegel upper half-space consisting of symmetric $g\times g$ complex matrices with positive 
definite imaginary part. For $k,m\in\mathbb{N}_{>0}$, we let $J_{k,m}^{\mathrm{cusp}}(\Gamma_{0})$ denote the complex vector 
space of Siegel--Jacobi cusp forms of weight $k$ and index $m$ for $\Gamma_{0}$. For any Siegel--Jacobi cusp form $\phi\in J_{k,
m}^{\mathrm{cusp}}(\Gamma_{0})$, its pointwise Petersson norm is defined for $(\tau,z)\in\mathbb{H}_{g}\times\mathbb{C}^{g}$ by 
the formula
\begin{displaymath}
\Vert\phi(\tau,z)\Vert_{\mathrm{Pet}}^{2}:=\vert\phi(\tau,z)\vert^{2}\,\det(\mathrm{Im}(\tau))^{k}\,e^{-4\pi m\mathrm{Im}(z)^{t}\mathrm
{Im}(\tau)^{-1}\mathrm{Im}(z)},
\end{displaymath}
which defines a $\Gamma_{0}\ltimes\mathbb{Z}^{2g}$-invariant function on $\mathbb{H}_{g}\times\mathbb{C}^{g}$. Assuming that 
$\phi$ is $L^{2}$-normalized with respect to the Petersson inner product on the quotient space $\Gamma_{0}\ltimes\mathbb{Z}^{2g}
\backslash\mathbb{H}_{g}\times\mathbb{C}^{g}$, we establish the bound
\begin{displaymath}
\Vert\phi\Vert_{L^{\infty}}:=\sup_{(\tau,z)\in\mathbb{H}_{g}\times\mathbb{C}^{g}}\Vert\phi(\tau,z)\Vert_{\mathrm{Pet}}=O_{\Gamma_
{0}, \epsilon}\big(k^{(3g^{2}+5g)/8}\,m^{g^{2}+5g/4+\epsilon}\big)
\end{displaymath}
for any $k\in\mathbb{Z}_{\geq g+1}$, $m\in\mathbb{Z}_{\geq 1}$, and a given $\epsilon > 0$, where the implied constant depends 
only on $\Gamma_{0}$ and the choice of $\epsilon>0$. 

The fundamental strategy of proof of this paper structurally mirrors the geometric and analytic approach established in the case $g=1$ 
in~\cite{AKP25}, showing that the core methodology carries over to higher dimensions once the proper symplectic machinery is deployed. 
The conceptual basis that allows us to execute the proof in higher genus is provided by the uniform bounds on average for Siegel cusp 
forms established by J.~Kramer and A.~Mandal~\cite{k1} via the associated Bergman kernels. By utilizing the estimates from~\cite{k1}, 
the problem is reduced to the study of the components of the Eichler--Zagier theta decomposition generalized to higher genus. More
specifically, we decompose the Siegel--Jacobi form into a finite linear combination of Siegel modular forms of half-integral weight $k-
1/2$ multiplied by higher-dimensional Jacobi theta functions. We then apply the maximum principle for subharmonic functions to bound 
the cusp forms of half-integral weight near the boundary of the Siegel fundamental domain, while controlling the growth of the Jacobi theta
functions through higher-dimensional integral tests. A slight technical distinction arises when comparing the resulting exponents in the 
case $g=1$, is that for arbitrary genus $g\geq 1$, the Siegel modular forms of half-integral weight appearing in the theta decomposition are 
considered with respect to the principal congruence subgroup $\Gamma(4m)$ rather than the larger group $\Gamma_{0,1}(4m)$ featured 
in the classical $g=1$ setting. The resulting larger index of the underlying congruence subgroup in $\Gamma_{0}$ increases the exponent 
of the index $m$ in our bound, leading to a slightly modified bound when evaluated at $g=1$. Nevertheless, combining these bounds via 
the Cauchy--Schwarz inequality yields the polynomial dependence on $k$ and $m$ stated above.

\subsection{Outline of the paper}
The paper is organized as follows. In Section~2, we collect the necessary background material concerning the geometry of the Siegel 
upper half-space $\mathbb{H}_{g}$, properties of the Bergman kernel for Siegel cusp forms of weight $k$ for $\mathrm{Sp}_{2g}
(\mathbb{R})$, and fundamental properties of Siegel--Jacobi cusp forms of weight $k$ and index $m$ for $\mathrm{Sp}_{2g}(\mathbb
{R})$. Section~3 is entirely devoted to the proof of our main result, for which we need the technical proofs of Lemma~\ref{lem1} and 
Lemma~\ref{lem2}, culminating in the proof of our main theorem (Theorem~\ref{thm3}).

\section{Background Material}

\subsection{Siegel cusp forms}

For $g\geq 1$, let 
\begin{align*}
\mathbb{H}_{g}:=\big\lbrace\tau=\xi+i\eta\,\big|\,\xi,\eta\in\mathrm{Sym}_{g}(\mathbb{R}),\eta>0\big\rbrace
\end{align*}
denote the Siegel upper half-space, where $\mathrm{Sym}_{g}(\mathbb{R})$ denotes the set of symmetric $g\times g$ matrices 
with real entries. Furthermore, let $\mu_{\mathrm{hyp}}$ denote the $(1,1)$-form associated to the natural metric on $\mathbb{H}_
{g}$, and let $\mu_{\mathrm{hyp}}^{\mathrm{vol}}$ denote the associated volume form, which for any $\tau=\xi+i\eta\in\mathbb{H}_
{g}$ is given by the  formula
\begin{align*} 
\mu_{\mathrm{hyp}}^{\mathrm{vol}}(\tau):=\frac{\displaystyle\bigwedge_{1\leq l\leq j\leq g}\mathrm{d}\xi_{l,j}\wedge\mathrm{d}
\eta_{l,j} }{\mathrm{det}(\eta)^{g+1}}.
\end{align*}

Let $\mathrm{M}_{2g}(\mathbb{R})$ denote the set of $g\times g$ matrices with real entries, and let
\begin{align*}
\mathrm{Sp}_{2g}(\mathbb{R}):=\big\lbrace \gamma\in \mathrm{M}_{2g}(\mathbb{R})\big|\, \gamma^t J_{2g}\gamma=J_{2g}\big
\rbrace,
\end{align*}
where 
\begin{align*}
J_{2g}:=
\left(\begin{matrix}
 0&\mathrm{id}_{g}\\- \mathrm{id}_{g}&0
\end{matrix}\right),
\end{align*} 
and $\mathrm{id}_{g}$ denotes the identity matrix of dimension $g$. 

Let $\Gamma\subset \mathrm{Sp}_{2g}(\mathbb{R})$ be a cofinite subgroup, which acts on $\mathbb{H}_{g}$ via fractional linear transformations. Let $X_{\Gamma}:=\Gamma\backslash\mathbb{H}_{g}$ denote the quotient space, which admits the structure of a complex orbifold of dimension $g(g+1)/2$, and is of finite hyperbolic volume. The hyperbolic metric $\mu_{\mathrm{hyp}}$ descends to $X_{\Gamma}$, to define a metric, which is compatible with its complex structure. 

For any $k\in \mathbb{R}_{\geq 1}$, let $S_{k,\chi}(\Gamma)$ denote the complex vector space of Siegel cusp forms of weight $k$ and character $\chi$ for $\Gamma$. The complex vector space $S_{k,\chi}(\Gamma)$ is endowed with a pointwise norm $\|\cdot\|_{\mathrm{Pet}}$, known as the Petersson norm, which for any $f\in S_{k,\chi}(\Gamma)$, and $\tau=\xi+i\eta\in \mathbb{H}_{g}$, is given by the  formula
\begin{align*}
\|f(\tau)\|_{\mathrm{Pet}}^{2}:=|f(\tau)|^{2}\mathrm{det}(\eta)^{k}.
\end{align*}

Furthermore, $S_{k,\chi}(\Gamma)$ is equipped with an $L^2$-norm $\langle\cdot,\cdot\rangle_{\mathrm{Pet}}$, known as the Petersson inner product, which for any $f,g\in S_{k,\chi}(\Gamma)$, is given by the formula
\begin{align*}
\langle f, g\rangle_{\mathrm{Pet}}:=\int_{\mathcal{F}_{\Gamma}}f(\tau)\overline{g(\tau)}\mathrm{det}(\eta)^{k}\mu_{\mathrm{hyp}}^{\mathrm{vol}}(\tau),
\end{align*}
where $\mathcal{F}_{\Gamma}\subset \mathbb{H}_{g}$ denotes a fixed fundamental domain of the complex manifold $X_{\Gamma}$. 

Let $\lbrace f_{1},\ldots, f_{d_{k}}\rbrace$ denote an orthonormal basis of $S_{k,\chi}(\Gamma)$ with respect to the Petersson inner product $\langle\cdot,\cdot\rangle_{\mathrm{Pet}}$, where $d_{k}$ is the  dimension of $S_{k,\chi}(\Gamma)$, as a complex vector space. Let $B_{k,\chi}$ denote the Bergman kernel associated to $S_{k,\chi}(\Gamma)$, which for any $\tau, \tau^{\prime}\in \mathbb{H}_{g}$, is given by the formula
\begin{align*}
B_{k,\chi}(\tau,\tau^{\prime})=\sum_{j=1}^{d_{k}}f_{j}(\tau)\overline{f_{j}(\tau^{\prime})}.
\end{align*}
From the Riesz representation theorem, it follows that the definition of $B_{k,\chi}$ is independent of the choice of basis for $S_{k,\chi}(\Gamma)$. 

\subsection{Estimates for cusp forms}\label{subsec-1.2}

Let $\Gamma_0\subset\mathrm{Sp}_{2g}(\mathbb{R})$ be a cofinite subgroup, and let $\Gamma\subset \Gamma_{0}$ be a finite index subgroup. With notations as above, for any $\tau=\xi+i\eta\in\mathbb{H}_{g}$, let $0<\lambda_1(\eta)\leq \cdots\leq \lambda_g(\eta)$ denote the eigenvalues of $\eta$. For $\lambda_{0}>0$, put 
\begin{align*}
&\mathcal{F}_{\Gamma,\lambda_{0}}=\big\lbrace \tau=\xi+i\eta\in \mathcal{F}_{\Gamma} \big|\,\lambda_g(\eta)\leq \lambda_{0}\big\rbrace;\\
&\mathcal{F}_{\Gamma}^{\lambda_{0}}=\big\lbrace \tau=\xi+i\eta\in \mathcal{F}_{\Gamma} \big|\,\lambda_g(\eta)\geq \lambda_{0}\big\rbrace.
\end{align*}
  
Adapting arguments from the proofs of Theorems~4.5 and ~4.6 of \cite{k1} to the setting of $k\in \mathbb{R}_{\geq g+1}$, for any $f\in S_{k,\chi}(\Gamma)$, which is normalized with respect to the Petersson inner product $\langle \cdot,\cdot\rangle_{\mathrm{Pet}}$, we have the bound
\begin{align}\label{b1}
\|f\|_{L^{\infty}}^{2}=\sup_{\tau\in \partial \mathcal{F}_{\Gamma}^{\lambda_{0}}}\|f(\tau)\|_{\mathrm{Pet}}^{2}\leq \sup_{\tau\in \partial\mathcal{F}_{\Gamma}^{\lambda_{0}}}\|B_{k,\chi}(\tau,\tau)\|_{\mathrm{Pet}}^{2}=O_{\Gamma_{0}}\big(k^{3g(g+1)/4}\big),
\end{align}
where $\lambda_{0}=k/c_{\Gamma}$ with a positive constant $c_{\Gamma}$ depending only on $\Gamma$; furthermore, the implied constant in the above bound, depends only on $\Gamma_{0}$.

\subsection{Jacobi cusp forms}

From now on, let $\Gamma_{0}:=\mathrm{Sp}_{2g}(\mathbb{Z})$ denote the full modular group, and let  $\mathcal{F}_{\Gamma_{0}}$ denote the standard Siegel fundamental domain of $X_{\Gamma_{0}}=\Gamma_0\backslash \mathbb{H}_{g}$. 
  
For $k,m\in\mathbb{N}_{>0}$, let $J_{k,m}^{\mathrm{cusp}}(\Gamma_{0})$ denote the space of Jacobi cusp forms of weight $k$ and index $m$ for $\Gamma_{0}=\mathrm{Sp}_{2g}(\mathbb{Z})$. Any $\phi\in J_{k,m}^{\mathrm{cusp}}(\Gamma_{0})$ is a holomorphic function 
\begin{align*}
\phi:\mathbb{H}_{g}\times \mathbb{C}^{g}\longrightarrow \mathbb{C},
\end{align*}
which for any $\left[\gamma=\begin{pmatrix}A& B\\C&D \end{pmatrix},\begin{pmatrix}\lambda\\ \mu\end{pmatrix}\right]\in \Gamma_{0}\ltimes\mathbb{Z}^{2g}$ and $(\tau,z)\in\mathbb{H}_{g}\times\mathbb{C}^{g}$, satisfies the following transformation formula
\begin{align*}
&\phi\big((A\tau+B)\cdot(C\tau+D)^{-1},(C\tau+D)^{-1}(z+\tau\lambda+\mu)\big)\mathrm{det}(C\tau+D)^{-k}\times\\&\mathrm{exp}\big(2\pi i m\big( \lambda^{t}\tau\lambda+2\lambda^{t}z-(z+\tau\lambda+\mu)^{t} (C\tau+D)^{-1}(z+\tau\lambda+\mu)\big) \big)=\phi(\tau,z).
\end{align*}

Furthermore, any $\phi\in J_{k,m}^{\mathrm{cusp}}(\Gamma_{0})$ admits the following Fourier expansion
\begin{align*}
\phi(\tau,z)=\sum_{\substack{N\in \mathrm{Sym}_{g}(\mathbb{Q}), N\text{ half-integral}\\  r\in\mathbb{Z}^{g},\,4mN-r\cdot r^{t}>0}} c(N,r)e^{2\pi i \mathrm{Tr}(N\tau)}\cdot e^{2\pi i r^{t}z}.
\end{align*}

For $\phi\in J_{k,m}^{\mathrm{cusp}}(\Gamma_{0})$, we define
\begin{align*}
\|\phi(\tau,z)\|_{\mathrm{Pet}}^{2}:=\big|\phi(\tau,z)\big|^{2}\mathrm{det}(\eta)^{k}e^{-4\pi m y^{t}\eta^{-1}y}\quad (\tau=\xi+i\eta\in\mathbb{H}_{g}, z=x+iy\in \mathbb{C}^{g}),
\end{align*}
which defines a $\Gamma_{0}\ltimes\mathbb{Z}^{2g}$-invariant function on $\mathbb{H}_{g}\times\mathbb{C}^{g}$, called the pointwise Petersson norm of $\phi$.
 
Let $\mathcal{D}_{\Gamma_{0}}$ denote a fixed fundamental domain of the quotient space $Y_{\Gamma_{0}}:= \Gamma_{0}\ltimes\mathbb{Z}^{2g}\backslash \mathbb{H}_{g}\times\mathbb{C}^{g}$, which admits the structure of a $g(g+3)/2$-dimensional complex manifold. The space $J_{k,m}^{\mathrm{cusp}}(\Gamma_{0})$ is equipped with the Petersson inner product 
\begin{align}\label{pet-ip}
\langle \phi,\psi\rangle_{\mathrm{Pet}}:=\int_{\mathcal{D}_{\Gamma_{0}}}\phi_{1}(\tau,z)\overline{\phi_{2}(\tau,z)}e^{-4\pi m y^{t}\eta^{-1}y}\frac{\displaystyle\bigwedge_{1\leq l\leq j\leq g}\mathrm{d}\xi_{l,j}\wedge\mathrm{d}\eta_{l,j} \displaystyle\bigwedge_{1\leq n \leq g} \mathrm{d}x_{n}\wedge\mathrm{d}y_{n}}{\mathrm{det}(\eta)^{2g+1-k}}.
\end{align}

\subsection{Theta functions and Eichler--Zagier relation}
 
A Jacobi form $\phi\in J_{k,m}^{\mathrm{cusp}}(\Gamma_{0})$ of weight $k$ and index $m$ for the Siegel modular group $\Gamma_{0}=\mathrm{Sp}_{2g}(\mathbb{Z})$ has an expansion of the form
\begin{align}\label{eich-zag}
\phi(\tau,z)=\sum_{\mu\in(\mathbb{Z}/2m\mathbb{Z})^{g}}h_{\mu}(\tau)\,\theta_{\mu,m}
(\tau,z),
\end{align}
where the $h_{\mu}$'s are Siegel modular forms of weight $k-1/2$ for some congruence subgroup $\Gamma$ containing the principal congruence subgroup $\Gamma(4m)$  of level $4m$, and 
the $\theta_{\mu,m}$'s are higher dimensional Jacobi theta functions given by the formula
\begin{align}\label{theta-defn}
\theta_{\mu,m}(\tau,z)=\sum_{n\in\mathbb{Z}^{g}}e^{2\pi im(n-\frac{\mu}{2m})^{t}\tau(n-\frac{\mu}{2m})+2\pi i (2mn-\mu)^{t}z}.
\end{align}

The Petersson norm of the theta function is given by the formula
\begin{align}\label{theta-pet-norm}
\|\theta_{\mu,m}(\tau,z)\|_{\mathrm{Pet}}:=|\theta_{\mu,m}(\tau,z)|\mathrm{det}(\eta)^{1/4}e^{-2\pi m y^{t}\eta^{-1}y}.
\end{align}
In fact, it is shown in Theorem 3.3 of \cite{ziegler} that decomposition \eqref{eich-zag} gives rise to the isomorphism
\begin{align*}
J_{k,m}^{\mathrm{cusp}}(\Gamma_{0})\cong V_{k-1/2} (\Gamma_0),
\end{align*}
where $V_{k-1/2} (\Gamma_0)$ denotes the complex vector space of vector-valued cusp forms of weight $k-1/2 $ with suitable transformation behaviour with respect to $\Gamma_0$.

For any 
\begin{align*}
\phi_{1}(\tau,z)=\sum_{\mu\in(\mathbb{Z}/2m\mathbb{Z})^{g}}h_{\mu,1}(\tau)\,\theta_{\mu,m}(\tau,z)\quad\mathrm{and}\quad \phi_{2}(\tau,z)=\sum_{\mu\in(\mathbb{Z}/2m\mathbb{Z})^{g}}h_{\mu,2}(\tau)\,\theta_{\mu,m}(\tau,z)
\end{align*}
be two Jacobi cusp forms of weight $k$ and index $m$ for $\Gamma_{0}$. Then, from decomposition \eqref{eich-zag}, 
we have the equality 
\begin{align}\label{dec-pet-ip}
\langle \phi_1,\phi_2\rangle_{\mathrm{Pet}}=\frac{1}{(4m)^{g/2}}\int_{\mathcal{F}_{\Gamma_{0}}}\sum_{\mu\in(\mathbb{Z}/2m\mathbb{Z})^{g}}h_{\mu,1}(\tau)\overline{h_{\mu,2}(\tau)}\mathrm{det}(\eta)^{k-1/2}\mu_{\mathrm{hyp}}^{\mathrm{vol}}(\tau),
\end{align}
see Lemma 3.4 in~\cite{ziegler}.

\section{Sup-norm bounds for Jacobi cusp forms}

With notations as above, for $k\in\mathbb{Z}_{\geq g+1}$, let $\phi\in J_{k,m}^{\mathrm{cusp}}(\Gamma_{0})$ be an $L^{2}$-normalized Jacobi cusp form of weight $k$ and index $m$ for the Siegel modular group $\Gamma_0=\mathrm{Sp}_{2g}(\mathbb{Z})$. 

We will bound the quantity 
\begin{align*}
\|\phi\|_{L^{\infty}}=\sup_{(\tau,z)\in\mathbb{H}_{g}\times \mathbb{C}^{g}}\|\phi(\tau,z)\|_{\mathrm{Pet}}.
\end{align*}

From the $\Gamma_{0}\ltimes\mathbb{Z}^{2g}$-invariance of the function $\|\phi(\tau,z)\|_{\mathrm{Pet}}$ on $\mathbb{H}_{g}\times \mathbb{C}^{g}$, it suffices to bound $\|\phi(\tau,z)\|_{\mathrm{Pet}}$, when $\tau$ ranges through a fundamental domain of $\mathcal{F}_{\Gamma_{0}}$ for $\Gamma_{0}$, and $z$ ranges through a fundamental domain of the corresponding abelian varieties 
$\mathbb{C}^{g}/(\mathbb{Z}^{g}\oplus \tau\mathbb{Z}^{g})$. 

Without loss of generality, we assume in the sequel that $\mathcal{F}_{\Gamma_0}\subset \mathbb{H}_{g}$, is a fixed fundamental domain of $\Gamma_0$, and is determined by the following properties: 
\begin{enumerate}
\item[(i)]
$\tau=\xi+i\eta\in\mathbb{H}_g$ with $\eta$ being Minkowski reduced. 
\item[(ii)]
$|\mathrm{det}(C\tau+D)|\geq 1$, for all $\begin{pmatrix}A & B \\ C &D\end{pmatrix}\in \Gamma_{0}$.
\item[(iii)]
For all $1\leq j,k\leq g$, the matrix $\xi=(\xi_{j,k})_{1\leq j,k\leq g}$ satisfies $|\xi_{j,k}|\leq 1/2$.
\item[(iv)]
For any $\eta=(\eta_{j,k})_{1\leq j,k\leq g}$, we have $\sqrt{3}/2\leq \eta_{1,1}\leq \cdots\leq \eta_{g,g}$.
\item[(v)]
There exists a constant $c(g)>0$ depending only on $g$, such that $\eta\geq c(g)\mathrm{id}_{g}$. 
\end{enumerate}

Furthermore, the fundamental domain of the abelian variety associated to a given $\tau\in \mathcal{F}_{\Gamma_0}$, is given by the following formula
\begin{align*}
A_{\tau}:=\big\lbrace z=x+iy\in\mathbb{C}^{g}\big|\,0\leq x_j\leq 1,\,0\leq y_{j}\leq \sum_{k=1}^{g}\eta_{j,k} ,\,\mathrm{for\,\,all\,\,}1\leq j\leq g\big\rbrace.
\end{align*}

To obtain our bound $\|\phi\|_{L^{\infty}}$, we start from the decomposition \eqref{eich-zag}
\begin{align*}
\phi(\tau,z)=\sum_{\mu\in(\mathbb{Z}/2m\mathbb{Z})^{g}}h_{\mu}(\tau)\,\theta_{\mu,m}(\tau,z).
\end{align*}

\begin{lem}\label{lem1}
With notations as above, for any $(\tau,z)\in\mathbb{H}_{g}\times \mathbb{C}^{g}$, we have the following bound
\begin{align*}
\sum_{\mu\in(\mathbb{Z}/2m\mathbb{Z})^{g}}\|\theta_{\mu,m}(\tau,z)\|_{\mathrm{Pet}}^{2}\leq (2m)^{g}\mathrm{det}(\eta)^{1/2} +O\big(m^{g-1/2}\,\mathrm{det}(\eta)^{1/2}\big).
\end{align*}
\begin{proof}
By definition, we have the equalities
\begin{align*}
\theta_{\mu,m}(\tau,z)&=\sum_{n\in\mathbb{Z}^{g}}e^{2\pi im(n-\frac{\mu}{2m})^{t}
\tau(n-\frac{\mu}{2m})+2\pi i(2mn-\mu)^{t}z} \\
&=\sum_{n\in\mathbb{Z}^{g}}e^{2\pi im(n^{t}\tau n)-2\pi i (n^{t}\tau\mu)+\frac{\pi i}{2m} (\mu^{t}\tau\mu)+4\pi im( n^{t}z)-2\pi i(\mu^{t} z)} \\
&=e^{\frac{\pi i}{2m}(\mu^{t}\tau\mu)-2\pi i(\mu^{t} z)}\,\sum_{n\in\mathbb{Z}^{g}}e^{2\pi im(n^{t}\tau n)+2\pi i(2m n^{t})(z-\frac{\tau\mu}{2m})} \\
&=\theta_{0,m}\Big(\tau,z-\frac{\tau\mu}{2m}\Big)\,e^{\frac{\pi i}{2m}(\mu^{t}\tau\mu)-2\pi i(\mu^{t} z)}.
\end{align*}
Taking absolute values, we derive
\begin{align*}
\vert\theta_{\mu,m}(\tau,z)\vert=\Big\vert\theta_{0,m}\Big(\tau,z-\frac{\tau\mu}{2m}\Big)\Big\vert\,e^{-\frac{\pi}{2m}(\mu^{t}\eta\mu)+2\pi(\mu^{t} y)}.
\end{align*}
This implies the relation
\begin{align*}
\Vert\theta_{\mu,m}(\tau,z)\Vert_{\mathrm{Pet}}&=\vert\theta_{\mu,m}(\tau,z)\vert\,\mathrm{det}(\eta)^{1/4}\,e^{-2\pi my^{t}\eta^{-1}y} \\[1mm]
&=\Big\vert\theta_{0,m}\Big(\tau,z-\frac{\tau\mu}{2m}\Big)\Big\vert\,\mathrm{det}(\eta)^{1/4}\,e^{-\frac{\pi}{2m}\mu^{t}\eta\mu+2\pi\mu^{t} y-2\pi m y^{t}\eta^{-1}y}
 \\[1mm]
&=\Big\vert\theta_{0,m}\Big(\tau,z-\frac{\tau\mu}{2m}\Big)\Big\vert\,\mathrm{det}(\eta)^{1/4}\,e^{-2\pi m(y-\frac{\eta\mu}{2m})^{t}\eta^{-1}(y-\frac{\eta\mu}{2m})}.
\end{align*}
From this, we derive for $\tau\in\mathcal{F}_{\Gamma_{0}}$ and $z\in A_{\tau}$, the bound
\begin{align}\label{lem1-eqn1}
\Vert\theta_{\mu,m}(\tau,z)\Vert_{\mathrm{Pet}}&=\Big\vert\theta_{0,m}\Big(\tau,z-\frac{\tau\mu}{2m}\Big)\Big\vert\,e^{-2\pi m(y-\frac{\eta\mu}{2m})^{t}\eta^{-1}(y-\frac{\eta\mu}{2m})}\,
\mathrm{det}(\eta)^{1/4}\notag\\[1mm]
&\le\sum_{n\in\mathbb{Z}^{g}}e^{-2\pi mn^{t}\eta n-4\pi mn^{t}(y-\frac{\eta\mu}{2m})-2\pi m(y-\frac{\eta\mu}{2m})\eta^{-1}(y-\frac{\eta\mu}{2m})}\,\mathrm{det}(\eta)^{1/4}\notag \\
&=\sum_{n\in\mathbb{Z}^{g}}e^{-2\pi m(\eta n+y-\frac{\eta\mu}{2m})^{t}\eta^{-1}(\eta n +y-\frac{\eta\mu}{2m})}\,\mathrm{det}(\eta)^{1/4}\notag\\
&=\sum_{n\in\mathbb{Z}^{g}}e^{-2\pi m(n+\eta^{-1}y-\frac{\mu}{2m})^{t}\eta(n +\eta^{-1}y-\frac{\mu}{2m})}\,\mathrm{det}(\eta)^{1/4}.
\end{align}
Since $\eta$ is a positive definite symmetric matrix, there exists an orthogonal matrix $P$ such that $\eta= P^{t}\mathrm{diag}
(\lambda_1,\cdots,\lambda_g)P$. Set 
\begin{align*}
\Lambda:=\bigg\lbrace n^{\prime}:= P\bigg(n+\eta^{-1}y-\frac{\mu}{2m}\bigg)\bigg|n\in\mathbb{Z}^{g}\bigg\rbrace,
\end{align*}
which spans a lattice in $\mathbb{R}^{g}$. Using inequality \eqref{lem1-eqn1}, and an integral test, we arrive for $\tau\in\mathcal{F}_{\Gamma_{0}}$ and $z\in A_{\tau}$ at the bound
\begin{align*}
&\Vert\theta_{\mu,m}(\tau,z)\Vert_{\mathrm{Pet}}\leq \sum_{n^{\prime}\in\Lambda}e^{-2\pi m n'^{t} \mathrm{diag}(\lambda_{1},\ldots,\lambda_{g})n^{\prime}}\,\mathrm{det}(\eta)^{1/4}
\notag\\[1mm]&\leq\prod_{j=1}^{g}\bigg(1+\int_{-\infty}^{\infty}e^{-2\pi m \lambda_{j}n_{j}'^{2}}\,dn'_{j}\bigg)\mathrm{det}(\eta)^{1/4} =\prod_{j=1}^{g}\bigg(1+\frac{1}{\sqrt{2m\lambda_j}}\bigg)\mathrm{det}(\eta)^{1/4},
\end{align*}
which leads to
\begin{displaymath}
\sum_{\mu\in(\mathbb{Z}/2m\mathbb{Z})^{g}}\Vert\theta_{\mu,m}(\tau,z)\Vert^{2}_{\mathrm{Pet}}\le (2m)^{g}\,\mathrm{det}(\eta)^{1/2}+O\big(m^{g-1/2}\,\mathrm{det}(\eta)^{1/2}\big).
\end{displaymath}
This completes the proof of the lemma.
\end{proof}
\end{lem}

\begin{lem}\label{lem2}
With notations as above, for any $k\in\mathbb{R}_{>g+1}$, let $f\in S_{k,\chi}(\Gamma(4m))$ be an $L^{2}$-normalized cusp form of weight $k$ and character $\chi$ for $\Gamma$. Then, we have the bound
\begin{align*}
\sup_{\tau\in\mathcal{F}_{\Gamma_{0}}}\,\|f(\tau)\|_{\mathrm{Pet}}^{2}\,\mathrm{det}(\eta)^{1/2}=O_{\Gamma_{0}}(k^{(3g^2+5g)/4}),
\end{align*}
where the implied constant depends on $\Gamma_0$.
\begin{proof}
For a given $L^{2}$-normalized $f\in S_{k,\chi}(\Gamma)$ and $\tau=\xi+i\eta\in \mathcal{F}_{\Gamma_{0}}$, we write
\begin{align*}
\Vert f(\tau)\Vert_{\mathrm{Pet}}^{2}\,\mathrm{det}(\eta)^{1/2}=
\bigg\vert\frac{f(\tau)}{e^{2\pi ic_{\Gamma(4m)}\mathrm{tr}(\tau)}}\bigg\vert^{2}\cdot\bigg\vert\frac{\mathrm{det}(\eta)^{k+1/2}}{
e^{2\pi  c_{\Gamma(4m)}\mathrm{tr}(\eta)}}\bigg\vert,
\end{align*}
where $c_{\Gamma(4m)}$ is as described in Subsection~\ref{subsec-1.2}. Now, the function $\big\vert f(\tau)/e^{2\pi ic_{\Gamma(4m)}\mathrm{tr}(\tau)}\big\vert^{2}$ is subharmonic and bounded in the neighborhood $\mathcal{F}_{\Gamma_{0}}^{\lambda_{0}}$, and thus takes its maximum on the boundary $\partial\mathcal{F}_{\Gamma_0}^{\lambda_{0}}$ by the maximum principle for subharmonic functions. 

Furthermore, from arguments from Theorem~4.5 from \cite{k1}, we find that the function 
\begin{align*}
\bigg\vert\frac{\mathrm{det}(\eta)^{k+1/2}}{e^{ 2\pi c_{\Gamma(4m)}\mathrm{tr}(\eta)}}\bigg\vert=\prod_{j=1}^{g}\frac{\lambda_{j}(\eta)^{k+1/2}}{e^{2\pi c_{\Gamma(4m)}\lambda_{j}(\eta)}}\qquad (\lambda_{j}(\eta)\text{ eigenvalues of }\eta) 
\end{align*} 
attains its maximum at 
\begin{align*}
\lambda_{1}(\eta)=\cdots=\lambda_{g}(\eta)=\frac{k+1/2}{c_{\Gamma(4m)}}.
\end{align*}
Now, for each $1\leq j \leq g$, the component $\lambda_{j}(\eta)^{k+1/2}/e^{2\pi c_{\Gamma(4m)}\lambda_{j}(\eta)}$ is a strictly monotonically decreasing function for $\lambda_{j}(\eta)>(k+1/2)/c_{\Gamma(4m)}$. By assuming for the moment that $k$ is large enough and by choosing
\begin{align*}
\lambda_{0}=\frac{k+1/2}{c_{\Gamma(4m)}}>\sqrt{3}/2,
\end{align*}
we derive from the arguments used in the proofs of Theorems~4.5 and~4.6 from \cite{k1}, in particular from bound~(4.24) that
\begin{align}
\notag
&\sup_{\tau\in\mathcal{F}_{\Gamma_{0}}}\big(\Vert f(\tau)\Vert_{\mathrm{Pet}}^{2}\,\mathrm{det}(\eta)^{1/2}\big)=\sup_{\tau\in\partial\mathcal{F}_{\Gamma_{0}}^{\lambda_{0}}}\big(\Vert f(\tau)\Vert_{\mathrm{Pet}}^{2}\,\mathrm{det}(\eta)^{1/2}\big) \\[2mm]
\notag
&\le\sup_{\tau\in\partial\mathcal{F}_{\Gamma_{0}}^{\lambda_{0}}}\big(\Vert B_{k,\chi}(\tau,\tau)\Vert_{\mathrm{Pet}}\,\mathrm{det}(\eta)^{1/2}\big)=O_{\Gamma_{0}}\big(k^{3g(g+1)/4}\cdot k^{g/2}\big)=O_{\Gamma_{0}}(k^{(3g^2+5g)/4}).
\end{align}
The remaining cases $k\in[g+1,\sqrt{3}/2\cdot c_{\Gamma(4m)}-1/2]$ can now be obviously added to the statement by possibly enlarging the implied constant.
\end{proof}
\end{lem}

\begin{thm}\label{thm3}
For $k\in\mathbb{Z}_{\geq g+1}$ and $m\in\mathbb{Z}_{\ge 1}$, let $\phi\in J_{k,m}^{\mathrm{cusp}}(\Gamma_{0})$ be an $L^{2}$-normalized Jacobi cusp form of weight $k$ and index $m$ for the full modular group $\Gamma_{0}=\mathrm{Sp}_{2g}(\mathbb{Z})$. Then, we have the $L^{\infty}$-norm 
bound
\begin{align*}
\Vert\phi\Vert^{2}_{L^{\infty}}=\sup_{(\tau,z)\in\mathbb{H}_{g}\times\mathbb{C}^{g}}\Vert\phi(\tau,z)\Vert^{2}_{\mathrm{Pet}}=O_{\Gamma_{0}, \epsilon}\big
(k^{(3g^{2}+5g)/4}\,m^{2g^{2}+5g/2+\epsilon}\big),
\end{align*}
where the implied constant depends on $\Gamma_{0}$ and the choice of $\epsilon>0$.
\begin{proof}
Substituting the decomposition~\eqref{eich-zag} of the Jacobi form $\phi(\tau,z)$ into its pointwise Petersson norm and applying the Cauchy--Schwarz inequality, we find the estimate
\begin{align*}
\notag
\Vert\phi(\tau,z)\Vert^{2}_{\mathrm{Pet}}&=\bigg\vert\sum_{\mu\in(\mathbb{Z}/2m\mathbb{Z})^{g}}h_{\mu}(\tau)\theta_{\mu,m}(\tau,z)\bigg\vert^{2}\,\mathrm{det}(\eta)^{k}\,e^{-4\pi my^{t}\eta^{-1}y} \\
\notag
&\leq\bigg(\sum_{\mu\in(\mathbb{Z}/2m\mathbb{Z})^{g}}\vert h_{\mu}(\tau)\vert^{2}\,\mathrm{det}(\eta)^{k-1/2}\bigg)\bigg(\sum_{\mu\in(\mathbb{Z}/2m\mathbb{Z})^{g}}\vert\theta_{\mu,m}(\tau,z)\vert^{2}\,\mathrm{det}(\eta)^{k-1/2}\,e^{-4\pi my^{t}\eta^{-1}y}\bigg) \\[1mm]
\notag
&=\bigg(\sum_{\mu\in(\mathbb{Z}/2m\mathbb{Z})^{g}}\Vert h_{\mu}(\tau)\Vert^{2}_{\mathrm{Pet}}\bigg)\bigg(\sum_{\mu\in(\mathbb{Z}/2m\mathbb{Z})^{g}}\Vert\theta_{\mu,m}(\tau,z)\Vert_{\mathrm{Pet}}^{2}\bigg).
\end{align*}
Thus, we arrive by means of Lemma~\ref{lem1} at the bound
\begin{align*}
\notag
\Vert\phi(\tau,z)\Vert^{2}_{\mathrm{Pet}}&\le\sum_{\mu\in (\mathbb{Z}/2m\mathbb{Z})^{g}}\Vert h_{\mu}(\tau)\Vert^{2}_{\mathrm{Pet}}\cdot\sum_{\mu\in(\mathbb{Z}/2m\mathbb{Z})^{g}}\Vert\theta_{\mu,m}(\tau,z)\Vert_{\mathrm{Pet}}^{2} \\
&\le (2m)^{g}\sum_{\mu\in (\mathbb{Z}/2m\mathbb{Z})^{g}}\Vert h_{\mu}(\tau)\Vert^{2}_{\mathrm{Pet}}\,\mathrm{det}(\eta)^{1/2}+O\bigg(m^{g-1/2}\,\sum_{\mu\in (\mathbb{Z}/2m\mathbb{Z})^{g}}\Vert h_{\mu}(\tau)\Vert^{2}_{\mathrm{Pet}}\bigg),
\end{align*}
from which we derive, after recalling the comment at the beginning of this section,
\begin{align}
\notag
\Vert\phi\Vert^{2}_{L^{\infty}}&=\sup_{(\tau,z)\in\mathbb{H}_{g}\times\mathbb{C}^{g}}\Vert\phi(\tau,z)\Vert^{2}_{\mathrm{Pet}}=\sup_{\substack{\tau\in\mathcal
{F}_{\Gamma_{0}}\\z\in A_{\tau}}}\Vert\phi(\tau,z)\Vert^{2}_{\mathrm{Pet}} \\
\label{thm3-eqn1}
&\le (2m)^{g}
%\sum_{\mu=0}^{2m-1}
\sum_{\mu\in(\mathbb{Z}/2m\mathbb{Z})^{g}}
\sup_{\tau\in\mathcal{F}_{\Gamma_{0}}}\big(\Vert h_{\mu}(\tau)\Vert^{2}_{\mathrm{Pet}}\,\mathrm{det}(\eta)^{1/2}\big)+O\bigg
(m^{g-1/2}\sum_{\mu\in(\mathbb{Z}/2m\mathbb{Z})^{g}}\sup_{\tau\in\mathcal{F}_{\Gamma_{0}}}\Vert h_{\mu}(\tau)\Vert^{2}_{\mathrm{Pet}}\bigg).
\end{align}
In order to apply bound~\eqref{b1} and Lemma~\ref{lem2} to the two summands in~\eqref{thm3-eqn1}, respectively, we need 
to $L^{2}$-normalize the modular forms $h_{\mu}$ ($\mu\in(\mathbb{Z}/2m\mathbb{Z})^{2g}$) under consideration. To do so, we observe that formula~\eqref{dec-pet-ip} in
conjunction with the $L^{2}$-normalization of  $\phi$ gives
\begin{align}
\label{thm3-eqn2}
\Vert\phi\Vert^{2}_{L^{2}}=\frac{1}{(4m)^{g/2}\,[\Gamma_{0}:\Gamma(4m)]}\sum_{\mu\in(\mathbb{Z}/2m\mathbb{Z})^{g}}\Vert h_{\mu}\Vert^{2}_{L^{2}}=1.
\end{align}
Moreover, for any $\epsilon>0$, we have the  bound
\begin{align}\label{thm3-eqn3}
[\Gamma_{0}:\Gamma(4m)]=O_{\epsilon}(m^{2g^{2}+g+\epsilon}),
\end{align}
where the implied constant depends only on choice of $\epsilon>0$. Combining bounds~\eqref{thm3-eqn2} and ~\eqref{thm3-eqn3}, we infer that
\begin{align*}
\sum_{\mu\in(\mathbb{Z}/2m\mathbb{Z})^{g}}\Vert h_{\mu}\Vert^{2}_{L^{2}}=(4m)^{g/2}\,[\Gamma_{0}:\Gamma(4m)]=O_{\epsilon}\big(m^{2g^{2}+3g/2+\epsilon}\big).
\end{align*}
With regard to the first summand in~\eqref{thm3-eqn1}, we thus obtain by Lemma~\ref{lem2} the bound
\begin{align}
\notag
&(2m)^{g}\sum_{\mu\in(\mathbb{Z}/2m\mathbb{Z})^{g}}\sup_{\tau\in\mathcal{F}_{\Gamma_{0}}}\big(\Vert h_{\mu}(\tau)\Vert^{2}_{\mathrm{Pet}}\,\mathrm{det}(\eta)^{1/2}\big)=\\[1mm]\notag
&(2m)^{g}\sum_{\mu\in(\mathbb{Z}/2m\mathbb{Z})^{g}}\sup_{\tau\in\mathcal{F}_{\Gamma_{0}}}\bigg(\frac{\Vert h_{\mu}(\tau)\Vert^{2}_{\mathrm{Pet}}}{\Vert h_{\mu}\Vert^{2}_{L^{2}}}\,\mathrm{det}(\eta)^{1/2}\bigg)\Vert h_{\mu}\Vert^{2}_{L^{2}} \\[1mm]
\label{thm3-eqn4}
&=O_{\Gamma_{0}}\big(m^{g}k^{(3g^{2}+5g)/4}\big)\sum_{\mu\in(\mathbb{Z}/2m\mathbb{Z})^{g}}\Vert h_{\mu}\Vert^{2}_{L^{2}}=O_{\Gamma_{0},\epsilon}\big(k^{(3g^{2}+5g)/4}\,m^{2g^{2}+5g/2+\epsilon}\big),
\end{align}
with an implied constant depending on $\Gamma_{0}$ and the choice of $\epsilon>0$. With regard to the second summand in~\eqref{thm3-eqn1},
we derive from Lemma~\ref{lem2} the bound 
\begin{align}
\notag
&m^{g-1/2}\sum_{\mu\in(\mathbb{Z}/2m\mathbb{Z})^{g}}\sup_{\tau\in\mathcal{F}_{\Gamma_{0}}}\Vert h_{\mu}(\tau)\Vert^{2}_{\mathrm{Pet}}=m^{g-1/2}\sum_{\mu\in(\mathbb{Z}/2m\mathbb{Z})^{g}}\sup_{\tau\in\mathcal{F}_{\Gamma_{0}}}\bigg(\frac{\Vert h_{\mu}(\tau)\Vert^{2}_{\mathrm{Pet}}}{\Vert h_{\mu}\Vert^{2}_{L^{2}}}\bigg)\Vert 
h_{\mu}\Vert^{2}_{L^{2}}= \\
\label{thm3-eqn5}&O_{\Gamma_{0}}\big(m^{g-1/2}k^{(3g^{2}+5g)/4}\big)\sum_{\mu\in(\mathbb{Z}/2m\mathbb{Z})^{g}}\Vert h_{\mu}\Vert^{2}_{L^{2}}=O_{\Gamma_{0},\epsilon}\big(k^{(3g^{2}+5g)/4}\,m^{2g^{2}+5g/2-1/2+\epsilon}\big),
\end{align}
with an implied constant depending on $\Gamma_{0}$ and the choice of $\epsilon>0$. The bounds~\eqref{thm3-eqn4} and~\eqref{thm3-eqn5} complete the proof of the theorem.
\end{proof}
\end{thm}


\begin{thebibliography}{AMM16}
\bibitem[AD23]{AD23}
P.~Anamby and S.~Das,
\newblock {\emph{Jacobi forms, Saito--Kurokawa lifts, their pullbacks and sup-norms on average}},
\newblock {Res. Math. Sci. \textbf{10}, 14 (2023)}.

\bibitem[AKP25]{AKP25}
A.~Aryasomayajula, J.~Kramer, and A.-M. von Pippich,
\newblock {\emph{$L^\infty$-norm bounds for Jacobi cusp forms}},
\newblock {Res. Number Theory \textbf{11}, 53 (2025)}.

\bibitem[EZ85]{EZ85}
M.~Eichler and D.~Zagier,
\newblock {\emph{The Theory of Jacobi Forms}},
\newblock {Progress in Mathematics, Vol. 55, Birkh\"auser, Boston, 1985}.

\bibitem[IS95]{IS95}
H.~Iwaniec and P.~Sarnak,
\newblock {\emph{$L^\infty$ norms of Hecke cusp forms}},
\newblock {Ann. of Math. (2) 141, 301--320 (1995)}.

\bibitem[KM23]{k1} J.~Kramer and A.~Mandal,
\newblock{\emph{Uniform sup-norm bounds on average for Siegel cusp forms}}, 
\newblock{arXiv:2310.05334v1 (2023)}.

\bibitem[Zi89]{ziegler} C.~Ziegler, 
\newblock{\emph{Jacobi Forms of Higher Degree}}, 
\newblock{Abh. Math. Sem. Univ. Hamburg 59, 191-224 (1989)}.

\end{thebibliography}
\end{document}